\documentclass[preprint,12pt]{elsarticle}

\usepackage{latexsym,amsmath,amssymb,xspace,comment,graphicx, url, color}
\usepackage{comment}
\usepackage[numbers]{natbib}

\journal{}

\begin{document}

\begin{frontmatter}



\title{Conversion of Mersenne Twister to double-precision floating-point numbers}


\author[harase]{Shin Harase\corref{cor1}}
\ead{harase@fc.ritsumei.ac.jp}
\address[harase]{College of Science and Engineering, Ritsumeikan University, 1-1-1 Nojihigashi, Kusatsu, Shiga, 525-8577, Japan.}

\cortext[cor1]{Corresponding author}

\begin{abstract}
The 32-bit Mersenne Twister generator MT19937 is a widely used random number generator. 
To generate numbers with more than 32 bits in bit length, and particularly when  
converting into 53-bit double-precision floating-point numbers in $[0,1)$ in the IEEE 754 format, 
the typical implementation concatenates 
two successive 32-bit integers and divides them by a power of $2$. 
In this case, the 32-bit MT19937 is optimized in terms of its equidistribution properties 
(the so-called dimension of equidistribution with $v$-bit accuracy) 
under the assumption that one will mainly be using 32-bit output values, 
and hence the concatenation sometimes degrades the dimension of equidistribution 
compared with the simple use of 32-bit outputs. 
In this paper, we analyze such phenomena by investigating  
hidden $\mathbb{F}_2$-linear relations among the bits of high-dimensional outputs. 
Accordingly, we report that MT19937 with a specific lag set
fails several statistical tests, 
such as the overlapping collision test, matrix rank test, and Hamming independence test.
\end{abstract}


\end{frontmatter}



\newtheorem{theorem}{Theorem}
\newtheorem{lemma}[theorem]{Lemma}
\newtheorem{corollary}[theorem]{Corollary}
\newdefinition{definition}{Definition}
\newdefinition{remark}{Remark}
\newproof{proof}{Proof}
\newtheorem{proposition}[theorem]{Proposition}

\section{Introduction} \label{sec:intro}

Random number generators (RNGs) (or more precisely, pseudorandom number generators)
are an essential tool for scientific computing, 
especially for Monte Carlo methods. 
Modern CPUs and operating systems are moving from 32 to 64 bits, 
so it is important to design and analyze RNGs in this respect.
Currently, the 32-bit Mersenne Twister (MT) generator MT19937, 
which was released by Matsumoto and Nishimura \cite{MT19937} in the 1990s, 
is a widely used RNG. 
MT19937 has been adopted in many programming languages and softwares, 
such as the header {\tt <random>} in the C++11 (and the later C++14) standard template library (STL). 

In applications, we often convert unsigned integers into 
53-bit double-precision floating-point numbers in $[0,1)$ in the IEEE 754 format \cite{citeulike:4001400}. 
This format is represented as a 64-bit word 
consisting of a sign (the most significant bit), an exponent (the next 11 bits), 
and a significand (the remaining 52 bits).
To generate numbers with more than 32 bits in bit length, 
one can either use an RNG implemented with 32-bit integers and concatenate successive
32-bit blocks to generate each output, or use an RNG whose
recurrence is implemented directly with 64 bits. 
For example, in the header {\tt <random>} of the C++11 STL in GCC, 
the former is adopted for 32-bit MT19937 \cite[Chapter 26.5]{ISO:2012:III}. 
Let $\mathbb{F}_2:= \{ 0, 1\}$ be the two-element field.  
Let
\begin{eqnarray} \label{eqn:usual}
\mathbf{x}_0, \mathbf{x}_1, \mathbf{x}_2, \mathbf{x}_3, \mathbf{x}_4, \mathbf{x}_5, \mathbf{x}_6, \mathbf{x}_7,\ldots \in \mathbb{F}_2^{32}
\end{eqnarray} 
be successive 32-bit integers generated by MT19937. 
To obtain 53-bit double-precision floating-point numbers in $[0, 1)$ 
(by selecting the 32-bit generator {\tt mt19937} and 
the distribution {\tt uniform\_real\_distribution(0.0,1.0)}) in the header {\tt <random>}, 
the function generates 64-bit unsigned integers 
\begin{eqnarray} \label{eqn:concatenation}
(\mathbf{x}_{1}, \mathbf{x}_{0}), (\mathbf{x}_{3}, \mathbf{x}_{2}), (\mathbf{x}_{5}, \mathbf{x}_{4}),  (\mathbf{x}_{7}, \mathbf{x}_{6}), \ldots \in \mathbb{F}_2^{64}
\end{eqnarray}
by concatenating two successive 32-bit integers in this order 
and divides them by the maximum value $2^{64}$. 

Here, 
the {\it dimension of equidistribution with $v$-bit accuracy} (denoted by $k(v)$) is used as 
a quality criterion for the high-dimensional uniformity of RNGs \cite{LP2009,MatsumotoSHN06}. 
Harase and Kimoto \cite{Harase:2018:IME:3175005.3159444} recently developed 
64-bit MT-type RNGs that are completely
optimized for this criterion (called {\it maximally equidistributed}). 
In the case of MT19937, Harase and Kimoto pointed out that 
the concatenation \eqref{eqn:concatenation} degrades $k(v)$ 
compared with the simple use of the 32-bit output \eqref{eqn:usual}. 

The aim of this paper is 
to analyze such phenomena in more detail by means of the method of Harase \cite{MR3191969}. 
Roughly speaking, MT19937 has hidden low-weight $\mathbb{F}_2$-linear relations 
(i.e., linear dependency relations among the bits whose sum becomes 0 over $\mathbb{F}_2$) 
in high-dimensional output values. 
When we concatenate the output values, 
the $\mathbb{F}_2$-linear relations are folded, 
so they decrease $k(v)$ and cause a deviation. 
In addition, we examine the behavior of the least significant bits of MT19937, 
and report that the concatenation \eqref{eqn:concatenation} with a specific lag set 
corresponding to the above $\mathbb{F}_2$-linear relations is 
rejected or has suspect $p$-values for several statistical tests in TestU01 \cite{MR2404400}, 
such as the overlapping collision test, matrix rank test, and Hamming independence test, as new findings. 
This rejection occurs in the dimension lower than that of the usual output \eqref{eqn:usual}. 

The remainder of this paper is organized as follows. 
In Section~\ref{sec:preliminaries}, we briefly recall the notation of the Mersenne Twister and the dimension of equidistribution with $v$-bit accuracy. 
In Section~\ref{sec:analysis}, 
we analyze the concatenation \eqref{eqn:concatenation} 
in terms of the dimension of equidistribution and $\mathbb{F}_2$-linear relations.  
In Section~4, we conduct statistical tests, and compare the output sequences \eqref{eqn:usual} and \eqref{eqn:concatenation}. 
In Section~\ref{sec:comparison}, we compare MT19937 and other RNGs. 
Section~\ref{sec:conclusion} concludes this paper. 

\section{Preliminaries} \label{sec:preliminaries}
\subsection{Mersenne Twister MT19937} \label{subsec:MT}
Matsumoto and Nishimura \cite{MT19937} developed the following RNG with a Mersenne prime period, known as the {\it Mersenne Twister} (MT). 
We perform addition and multiplication over the two-element field $\mathbb{F}_2$ (or modulo 2). 
Let $w$ be the word size of the intended machine. 
For initial seeds $\mathbf{w}_0, \mathbf{w}_1, \ldots, \mathbf{w}_{N-1} \in \mathbb{F}_2^w$, 
we generate a $w$-bit vector sequence $\mathbf{w}_N, \mathbf{w}_{N+1}, \mathbf{w}_{N+2},$
$\ldots \in \mathbb{F}_2^w$ based on the following recursion:
\begin{eqnarray} \label{eqn:MT_recurtion}
\mathbf{w}_{i+N} := \mathbf{w}_{i+M} \oplus (\overline{\mathbf{w}}_i^{w-r} \mid \underline{\mathbf{w}}_{i+1}^{r}) A \quad (i = 0, 1, 2, \ldots ).
\end{eqnarray}
Here, $\oplus$ is the bitwise exclusive-OR (i.e., addition in $\mathbb{F}_2^w$), 
$\overline{\mathbf{w}}_i^{w-r} \in \mathbb{F}_2^{w-r}$ and $\underline{\mathbf{w}}_{i+1}^{r} \in \mathbb{F}_2^r$ denote the $w-r$ most significant bits (MSBs) of $\mathbf{w}_i$ and 
$r$ least significant bits (LSBs) of $\mathbf{w}_{i+1}$, respectively, 
and the $w$-bit vector 
$(\overline{\mathbf{w}}_i^{w-r} \mid \underline{\mathbf{w}}_{i+1}^{r}) $
 is the concatenation of  
the $(w-r)$-bit vector $\overline{\mathbf{w}}_i^{w-r}$ and the $r$-bit vector $\underline{\mathbf{w}}_{i+1}^{r}$ in that order.  
$N$ is an integer, $M$ is an integer such that $0 < M < N-1$, and $A \in \mathbb{F}_2^{w \times w}$ is a $w \times w$ regular matrix. 
Further, to improve the randomness (especially for the dimension of equidistribution $k(v)$, as described below), 
we consider a $w$-bit output sequence
\begin{eqnarray} \label{eqn:MT_tempering}
\mathbf{w}_{N}T, \mathbf{w}_{N+1}T, \mathbf{w}_{N+2}T, \mathbf{w}_{N+3}T, \ldots \in \mathbb{F}_2^w
\end{eqnarray}
by multiplying a suitable $w \times w$ regular matrix $T \in \mathbb{F}_2^{w \times w}$ by the sequence from (\ref{eqn:MT_recurtion}).  
Matsumoto and Nishimura \cite{MT19937} proposed the parameter values 
$(w, N, M, r) = (32, 624, 397, 31)$ and quickly computable matrices $A$ and $T$  
so as to attain the (maximal) period $2^p-1$ ($p := Nw-r = 19937$), which is a Mersenne prime. 
This algorithm is called the {\it Mersenne Twister} generator {\it MT19937}. 
We identify the output (\ref{eqn:MT_tempering}) as the 32-bit integer sequence \eqref{eqn:usual}. 

\subsection{Dimension of equidistribution with $v$-bit accuracy $k(v)$} \label{subsec:equidistribution} 

The {\it dimension of equidistribution with $v$-bit accuracy} is widely used 
as a criterion for $\mathbb{F}_2$-linear RNGs \cite{LP2009,MatsumotoSHN06}. 
Let ${\rm tr}_v: \mathbb{F}_2^w \rightarrow \mathbb{F}_2^v$ be a truncation function 
that takes the $v$ MSBs. 
\begin{definition}[Dimension of equidistribution with $v$-bit accuracy] \label{def:equidistribution}
A sequence of $w$-bit integers with period $P = 2^p-1$
\begin{eqnarray} \label{eqn:periodic}
\mathbf{x}_0, \mathbf{x}_1, \ldots, \mathbf{x}_{P-1}, \mathbf{x}_P = \mathbf{x}_0, \ldots
\end{eqnarray}
is said to be {\it $k$-dimensionally equidistributed with $v$-bit accuracy} if the successive 
$k$-tuples of the $v$ MSBs 
\begin{eqnarray} \label{eqn:kv bits}
 ({\rm tr}_v(\mathbf{x}_i), {\rm tr}_v(\mathbf{x}_{i+1}), \ldots, {\rm tr}_v(\mathbf{x}_{i+k-1})) \in \mathbb{F}_2^{vk} \qquad (i = 0, \ldots, P-1)
\end{eqnarray}
are uniformly distributed over all possible $kv$-bit patterns except for the all-$0$ pattern, 
which occurs once less often, that is, each distinct $k$-tuple of $v$-bit words appears 
the same number of times in the sequence (with the exception of the $k$-tuple of zeros). 
The largest value of $k$ with this property is called the {\it dimension of equidistribution 
with $v$-bit accuracy}, denoted by $k(v)$. 
\end{definition}
As a criterion of uniformity, larger values of $k(v)$ for each $1 \leq v \leq w$ are desirable \cite{Tootill}. 
For $P = 2^p-1$, there is an upper bound $k(v) \leq \lfloor p/v \rfloor$. 
The gap $d(v) := \lfloor p/v \rfloor -k(v) \geq 0$ is called the {\it dimension defect at $v$}, and the sum of the gaps 
$\Delta:=\sum_{v = 1}^{w} d(v)$ is called the {\it total dimension defect}. 
If $\Delta = 0$, the generator is said to be {\it maximally equidistributed}. 
For generators such as MT 
based on linear recurrences over $\mathbb{F}_2$, 
which are called {\it $\mathbb{F}_2$-linear generators} \cite{LP2009,MatsumotoSHN06}, 
one can quickly compute $k(v)$ for $v = 1, \ldots, w$ 
using lattice reduction algorithms over formal power series fields \cite{CL2000,MR2827396}. MT19937 has $\Delta = 6750$. 
The sequence \eqref{eqn:concatenation} can also be viewed as 
a 64-bit $\mathbb{F}_2$-linear generator, 
so $k(v)$ is computable, see \cite[Remark~4.2]{Harase:2018:IME:3175005.3159444} for details.

\section{Concatenation of output values of MT19937} \label{sec:analysis}
We compute $k(v)$ for 32- and 64-bit output sequences \eqref{eqn:usual} and \eqref{eqn:concatenation}.  
Table~\ref{table:MT} lists $k(v)$ for $v = 1, \ldots, 32$. 
For the 64-bit sequence \eqref{eqn:concatenation}, 
we have $k(33) = \cdots = k(48) = 312, k(49) = \cdots = k(64) = 311$, and $\Delta = 13527$. 
In each case, $k(v)$ fluctuates slightly for $v = 1, \ldots, 11$. 
However, for $v = 12, \ldots, 16$, 
we have $k(v) = 623 < 1246$ in \eqref{eqn:concatenation}, 
and these values of $k(v)$ are half those of the 32-bit sequence \eqref{eqn:usual}. 
If the generator is maximally equidistributed, the optimal value of $k(12)$ is $\lfloor 19937/12 \rfloor = 1661$.

\begin{table}
\caption{Dimension of equidistribution $k(v)$ for 32-bit and 64-bit output sequences \eqref{eqn:usual} and \eqref{eqn:concatenation} of MT19937 $(v = 1, \ldots, 32)$.}
\label{table:MT}
{\footnotesize
\begin{tabular}{|c|r|r||c|r|r||c|r|r||c|r|r|} \hline
$v$ & 64-bit & 32-bit & $v$ & 64-bit & 32-bit & $v$ & 64-bit & 32-bit & $v$ & 64-bit & 32-bit \\ \hline
$1$ & $19937$ & $19937$  & $9$ & $2068$ & $1869$ & $17$ & $623$ & $623$ & $25$ & $623$ & $623$\\ \hline
$2$ & $9968$ & $9968$ & $10$ & $1869$ & $1869$ & $18$ & $623$ & $623$ & $26$ & $623$ & $623$ \\ \hline
$3$ & $ 6643$ & $6240$ & $11$ & $1558$ & $1248$ & $19$ & $623$ & $623$ & $27$ & $623$ & $623$ \\ \hline
$4$ & $4983$ & $4984$ & $12$ & $623$ & $1246$ & $20$ & $623$ & $623$ & $28$ & $623$ & $623$ \\ \hline
$5$ & $3894$ & $3738$ & $13$ & $623$ & $1246$ & $21$ & $623$ & $623$ & $29$ & $623$ & $623$ \\ \hline
$6$ & $2917$ & $3115$ & $14$ & $623$ & $1246$ & $22$ & $623$ & $623$ & $30$ & $510$ & $623$ \\ \hline
$7$ & $2294$ & $2493$ & $15$ & $623$ & $1246$ & $23$ & $623$ & $623$ & $31$ & $510$ & $623$ \\ \hline
$8$ & $2180$ & $2492$ & $16$ & $623$ & $1246$ & $24$ & $623$ & $623$ & $32$ & $510$ & $623$ \\ \hline
\end{tabular}}
\end{table}

We can explain this phenomenon as follows. 
Let us denote the components of $\mathbf{x}_i$ by $\mathbf{x}_i =: (x_{i,0}, \ldots, x_{i, 31})$. 
For the 32-bit output sequence \eqref{eqn:usual}, we construct a polynomial lattice associated with the $\mathbb{F}_2$-vector sequence, 
which was proposed by Couture and L'Ecuyer \cite{CL2000}, 
and analyze the lattice structure by means of the method of Harase \cite[Section~5]{MR3191969}. 
According to Eq.~(13) in \cite{MR3191969}, 
we have a five-term $\mathbb{F}_2$-linear relation on the 12 MSBs of \eqref{eqn:usual}:
\begin{eqnarray} \label{eqn:relation1}
x_{i,2} + x_{i+792,4} +x_{i+792,11} + x_{i+1246,4} +x_{i+1246,11} = 0.
\end{eqnarray}
When we concatenate the output values \eqref{eqn:usual} as in \eqref{eqn:concatenation}, 
the $\mathbb{F}_2$-linear relation \eqref{eqn:relation1} is folded into half of the dimension.
This destroys the uniformity in \eqref{eqn:kv bits},   
and hence $k(12)$ decreases in \eqref{eqn:concatenation}. 

Conversely, Haramoto et al.~\cite{MR2344296,MR3145575} 
recently pointed out that several generators fail statistical tests for the LSBs, 
and so research on RNGs has begun to focus on the behavior of the LSBs. 
Thus, we investigate $\mathbb{F}_2$-linear relations for the LSBs of MT19937. 
For this, we apply the method of Harase \cite{MR3191969} to the reverse generator, i.e., the generator obtained by reversing the order of the 32 bits in \eqref{eqn:usual}. 
As a result, we have $\Delta = 14850$ for the 32-bit reverse generator 
and detect a five-term $\mathbb{F}_2$-linear relation between the bits 21 to 30 (i.e., discarding the first 20 bits and taking the next 10) of any 1247 successive values in \eqref{eqn:usual}:
\begin{eqnarray} \label{eqn:relation2}
x_{i, 20} + x_{i+ 792, 22} + x_{i + 792, 29} + x_{i + 1246, 22} + x_{i + 1246, 29}  = 0.
\end{eqnarray}

When we have $\mathbb{F}_2$-linear relations 
with a small number of terms (i.e., the number of terms $\leq 6$), 
it may be possible to detect deviations in some statistical tests. 
Matsumoto and Nishimura \cite{MR1958868} gave a theoretical justification 
for this fact in terms of coding theory (see Remark~\ref{remark:coding theory}).
In Section~\ref{sec:testing}, 
we derive some deviations based on \eqref{eqn:relation1} and \eqref{eqn:relation2} among $\{ \mathbf{x}_{i}, \mathbf{x}_{i+792}, \mathbf{x}_{i+1246} \}$ for MT19937.
\begin{remark} \label{remark:linear_relation}
Harase \cite{MR3191969} pointed out that 
MT19937 has the six-term $\mathbb{F}_2$-linear relation
\begin{eqnarray} \label{eqn:relation3}
x_{i,1} +x_{i,16} +x_{i+396,2} +x_{i+396,17} + x_{i+623,2} + x_{i+623,17} =0
\end{eqnarray}
and several low-weight $\mathbb{F}_2$-linear relations among $\{ \mathbf{x}_i, \mathbf{x}_{i+396}, \mathbf{x}_{i+623} \}$. 
In the sequence \eqref{eqn:concatenation}, 
the middle term $\mathbf{x}_{i+396}$ is shifted to the 32 LSBs, 
so the $\mathbb{F}_2$-linear relation \eqref{eqn:relation3} 
no longer holds for the 32 MSBs of \eqref{eqn:concatenation}.
\end{remark}

\begin{remark}
One might consider the concatenation
\begin{eqnarray} \label{eqn:concatenation2}
(\mathbf{x}_{0}, \mathbf{x}_{1}), (\mathbf{x}_{2}, \mathbf{x}_{3}), (\mathbf{x}_{4}, \mathbf{x}_{5}),  (\mathbf{x}_{6}, \mathbf{x}_{7}), \ldots \in \mathbb{F}_2^{64},
\end{eqnarray}
which is arranged in the reverse order of \eqref{eqn:concatenation}. 
In this case, 
the $k(v)$ values for $v = 1, \ldots, 32$ are the same as \eqref{eqn:concatenation} because the 32 MSBs of \eqref{eqn:concatenation2} are 
from every other sampling of \eqref{eqn:usual}, but the $k(v)$ for $v = 32, \ldots, 64$ are slightly different: $k(33) = \cdots = k(64) = 311$ and $\Delta = 13543$, see \cite[Remark~4.2]{Harase:2018:IME:3175005.3159444}. 
This result implies that we cannot obtain a small $\Delta$ 
by the simple concatenation of the successive 32-bit integers of MT19937 to produce 64-bit integers.
\end{remark}

\begin{remark} \label{remark:mt19937ar}
Matsumoto and Nishimura later released the 2002 version of MT ({\tt mt19937ar}), which
has a better initialization procedure than the original 1998 version \cite{MT19937}. 
In this generator, 
the function ${\tt genrand\_res53()}$ is added 
to generate double-precision floating-point numbers in $[0,1)$ with 53-bit resolution. 
This function sets the 32-bit integers 
$a \gets (\mathbf{x}_{2i} \gg 5)$ and $b \gets (\mathbf{x}_{2i+1} \gg 6)$ 
and returns a real number $(a \times 67108864.0+b) \times (1.0/9007199254740992.0) \in [0,1)$ 
using the usual integer arithmetic, where $(\mathbf{x} \gg l)$ denotes a right (i.e., zero-padded) shift by $l$ bits. 
Namely, the function extracts the 27 MSBs of $\mathbf{x}_{2i}$ and the 26 MSBs of $\mathbf{x}_{2i+1}$, 
concatenates them, and converts the result into a 53-bit double-precision floating-point number in $[0,1)$ 
by dividing by $2^{53}$. 
In this case, the $k(v)$ values for $v = 1, \ldots, 26$ are the same as in 
\eqref{eqn:concatenation} and \eqref{eqn:concatenation2}, 
but $k(27) = \cdots = k(29) = 623$, $k(30) = 425$, and $k(31) = \cdots = k(52) = 311$. 
Note that this conversion is exactly the same as 
\eqref{eqn:concatenation} and \eqref{eqn:concatenation2} for the 26 MSBs, 
but it avoids the $\mathbb{F}_2$-linear relation \eqref{eqn:relation2} for the LSBs. 
\end{remark}

\begin{remark} \label{remark:coding theory}
In the case of $\mathbb{F}_2$-linear generators \cite{LP2009,MatsumotoSHN06}, including MT, 
there is no constant $\mathbb{F}_2$-linear relation among the $(v \times k(v))$ bits in \eqref{eqn:kv bits}, 
and it is expected that these bits have no deviation and pass statistical tests. 
If $k \geq k(v)+1$, we have an $\mathbb{F}_2$-linear relation among the $(v \times k)$ bits in \eqref{eqn:kv bits}, 
such as \eqref{eqn:relation1}, \eqref{eqn:relation2}, or \eqref{eqn:relation3}. 
According to \cite{MR1958868}, if we have an $\mathbb{F}_2$-linear relation 
with a small number of terms (i.e., the number of terms $\leq 6$), 
there is a possibility of a strong dependency among the bits in the RNG, leading to frequent failures in some statistical tests, 
but if the minimum number of terms is $\geq 20$, it is difficult to detect the dependency in modern computers because of bit mixing. 
Thus, we should avoid such a low-weight $\mathbb{F}_2$-linear relation. 
The problem for checking low-weight $\mathbb{F}_2$-linear relations is reduced to the weight enumeration of a linear code, 
which is generally NP-complete. 
Harase \cite{MR3191969} proposed an efficient method
for finding a basis of $\mathbb{F}_2$-linear relations (i.e., a basis of a linear code) using lattice basis reduction algorithms,
especially in the case of $k = k(v)+1$, i.e., just one dimension beyond $k(v)$. 
For more details, we refer the reader to \cite[Section~3]{MR1958868} and \cite[Section~5]{MR3191969}. 
\end{remark}

\section{Empirical statistical tests} \label{sec:testing}
L'Ecuyer and Touzin \cite{L'Ecuyer:2004:DRN:961292.961302} 
and L'Ecuyer and Simard \cite{MR3246605} reported that certain multiple recursive generators 
with specific lags have poor lattice structures and fail statistical tests, 
such as the birthday spacings test \cite{Knuth:1997:ACP:270146,MR1823110,Marsaglia1985}. 
Harase \cite{MR3191969} reported that MT19937 with specific lag sets 
also fails the birthday spacings test in a similar manner. 
In this section, we show that there exist other statistical tests 
in the TestU01 package \cite{MR2404400} for which the output \eqref{eqn:usual} passes but 
\eqref{eqn:concatenation} fails under a certain condition. 

The TestU01 package is inherently a 32-bit suite. Thus, for 32-bit MT19937, 
to generate  a sequence $u_0, u_1, u_2, u_3, \ldots \in [0, 1)$, 
we implement the following two conversions:
\begin{enumerate}
\item[(a)] Divide 32-bit integers $\mathbf{x}_i$ in \eqref{eqn:usual} by $2^{32}$;
\item[(b)] Divide 64-bit integers $(\mathbf{x}_{2i+1}, \mathbf{x}_{2i})$ in \eqref{eqn:concatenation} by $2^{64}$.
\end{enumerate} 
We further consider the sequence with a lag set $I = \{ j_1, \ldots, j_t \}$:
\begin{eqnarray} \label{eqn:lac_real}
u_{j_1}, \ldots, u_{j_t}, \ldots, u_{(j_t+1)i + j_1}, \ldots, u_{(j_t+1)i + j_t}, \ldots \in [0,1).
\end{eqnarray}
For some tests, each output value $u_i$ is replaced by $2^\tau u_{i} \mod 1$
to skip the first $\tau$ bits and test the LSBs. 
If $\tau= 0$, this is the usual output sequence. 

In the following, we test the null hypothesis $\mathcal{H}_0$: the RNG output \eqref{eqn:lac_real} 
is perfectly random.  
Throughout our tests, we set $I = \{ 0, 396, 623 \}$. 
Note that these indices correspond to \eqref{eqn:relation3}, 
and are half the values in \eqref{eqn:relation1} and \eqref{eqn:relation2}; 
this is because the 64-bit sequence \eqref{eqn:concatenation} is 
the concatenation of two successive 32-bit integers in \eqref{eqn:usual} and  
the $\mathbb{F}_2$-linear relations \eqref{eqn:relation1} and \eqref{eqn:relation2} are folded. 
We use the default test parameter sets in the SmallCrush, Crush, and BigCrush batteries of TestU01, 
and conduct tests for five different initial seeds. 
Note that MT19937 (or $\mathbb{F}_2$-linear generators in general)  
fails the {\it linear complexity} test 
and the {\it matrix rank} test in ado-hoc setting that uses large matrices 
because of $\mathbb{F}_2$-linear artifacts 
(see also Table~\ref{table:mt_matrixrank} and the last sentence before Remark~\ref{remark:decimation}).

Harase \cite{MR3191969} pointed out that, for Test 12 of Crush in TestU01, 
conversions (a) and (b) with $I = \{ 0, 396, 623 \}$ in \eqref{eqn:lac_real} 
are rejected or have suspect $p$-values for the {\it birthday spacings test} 
\cite{Knuth:1997:ACP:270146,Marsaglia1985,MR1823110}.
(More precisely, instead of conversion (b) with $I = \{ 0, 396, 623 \}$, 
he showed the deviation of conversion (a) with $I = \{ 0, 792, 1246 \}$, 
see Remark~\ref{remark:decimation} for details.) 
In TestU01, Test 14 of BigCrush is the birthday spacings test 
with 20 replications (i.e., the number of two-level tests) increasing the statistical power of Test 12 of Crush with 5 replications. 
For comparison, we examine the results. 
Let us select two positive numbers $n$, $t$ 
and generate the ``non-overlapping" $t$-dimensional points $\mathbf{u}_0, \ldots, \mathbf{u}_{n-1} \in [0,1)^t$, where $\mathbf{u}_i := (u_{it}, \dots, u_{it+t-1})$ for $i = 0, \ldots, n-1$. 
We partition the hypercube $[0,1)^t$ 
into $d^t$ cubic boxes of equal size by dividing $[0, 1)$ into $d$ equal segments.
Let $K_1 \leq K_2 \leq \cdots \leq K_n$ be the numbers of the boxes into which  
these points have fallen, sorted by increasing order. 
Define the spacings $S_j := K_{j+1} - K_j$ for $j = 1, \ldots, n - 1$. 
Let $Y$ be the total number of collisions of these spacings, i.e., 
the number of values of $j \in \{1, \cdots, n-2\}$ such that $S_{(j + 1)} = S_{(j)}$, 
where $S_{(1)}, \ldots, S_{(n-1)}$ are the spacings sorted by increasing order. 
If $d^t$ is large and $n^3/(4d^t)$ is not too large,
$Y$ is approximately a Poisson random variable with mean $n^3/(4d^t)$ under $\mathcal{H}_0$. 
We repeat the tests $\tilde{N}$ times, count the number $Y$ of collisions, and compute the (right) $p$-value 
using the sum. 
Table~\ref{table:mt_birthday} summarizes the $p$-values for 
$(\tilde{N}, n, \tau, d, t) = (20, 2 \times 10^7, 0, 2^{21}, 3)$ in Test 14 of BigCrush. 
In this test, the three-dimensional outputs with 21-bit accuracy are investigated, 
so both (a) and (b) with $I = \{ 0, 396, 623 \}$ are 
rejected in accordance with \eqref{eqn:relation3} and \eqref{eqn:relation1}, respectively. 

\begin{table}[t]
\centering
\caption{Summary of $p$-values on the birthday spacings test (Test 14 of BigCrush) for MT19937 with $I = \{ 0, 396, 623 \}$ }
\label{table:mt_birthday}
{\footnotesize
\begin{tabular}{l||c|c|c|c|c} \hline
{Conversion} & {1st} & {2nd} & {3rd} & {4th} & {5th} \\ \hline  \hline
{32-bit (a)} & $3.3 \times 10^{-52}$ & $1.1 \times 10^{-73}$ & $1.8 \times 10^{-57}$ & $7.3 \times 10^{-63}$ & $8.8 \times 10^{-44}$ \\ \hline 
{64-bit (b)} & $3.8 \times 10^{-19}$ & $6.0 \times 10^{-21}$ & $2.8 \times 10^{-17}$ & $6.8 \times 10^{-15}$ & $4.3 \times 10^{-19}$ \\ \hline 
\end{tabular}}
\end{table}

Next, we conduct the {\it overlapping collision} test \cite[Sections 2 and 3]{MR1951060}. 
Let us generate the ``overlapping" $t$-dimensional points $\mathbf{v}_0, \ldots, \mathbf{v}_{n-1} \in [0,1)^t$,  
where the points are defined by 
$\mathbf{v}_0 = (u_0, \ldots, u_{t-1}), \mathbf{v}_1 = (u_1, \ldots,  u_t), \ldots, \mathbf{v}_{n-2} = (u_{n-2}, u_{n-1}, u_0, \ldots, u_{t-3}), \mathbf{v}_{n-1} = (u_{n-1}, u_0, \ldots, u_{t-2})$. 
Consider the above hypercube and its partition. 
Let $C$ be the number of collisions, i.e., 
the number of times a point falls in a cell that already has one or more points. 
Under $\mathcal{H}_0$, $C$ is approximately a Poisson random variable
with mean $n^2/(2d^t)$ if $n/d^t \leq 1$. 
We repeat the tests $\tilde{N}$ times, count the number $C$ of collisions, and compute the $p$-value. 
We use the parameter sets 
$(\tilde{N}, n, \tau, d, t) = (30, 2 \times 10^7, 0, 2^{14}, 3)$ and $(30, 2 \times 10^7, 16, 2^{14}, 3)$ from Tests 5 and 6 of BigCrush, respectively. 
Table~\ref{table:mt_collision} summarizes the $p$-values. 
For $I = \{ 0, 396, 623 \}$, 
conversion (a) passes the tests, but conversion (b) fails Test 5  
and has suspect $p$-values for Test 6. 
These results agree with 
the nonexistence of 
$\mathbb{F}_2$-linear relations and 
the presence of the $\mathbb{F}_2$-linear relations in \eqref{eqn:relation1} and \eqref{eqn:relation2}. 

\begin{table}[t]
\centering
\caption{Summary of $p$-values on the overlapping collision tests (Tests 5 and 6 of BigCrush) for MT19937 with $I = \{ 0, 396, 623\}$} \label{table:mt_collision}
{\footnotesize
\begin{tabular}{l||c|c|c|c|c} \hline
{Conversion} & {1st} & {2nd} & {3rd} & {4th} & {5th} \\ \hline  \hline
{32-bit (a) (Test 5)} & $0.28$ & $0.58$ & $0.78$ & $0.27$ & $0.25$ \\ \hline 
{64-bit (b) (Test 5)} & $1.8 \times 10^{-25}$ & $5.7 \times 10^{-34}$ & $1.7 \times 10^{-23}$ & $2.8 \times 10^{-37}$ & $8.8 \times 10^{-39}$ \\ \hline
{32-bit (a) (Test 6)} & $0.32$ & $0.21$ & $0.51$ & $0.32$ & $0.32$ \\ \hline 
{64-bit (b) (Test 6)} & $1.4 \times 10^{-4}$ & $7.6 \times 10^{-7}$ & $2.8 \times 10^{-17}$ & $4.8 \times 10^{-8}$ & $1.1 \times 10^{-5}$ \\ \hline
\end{tabular}}
\end{table}

To detect a deviation of the LSBs, 
TestU01 has several tests that skip the $\tau = 20$ 
MSBs and look at the next $\sigma = 10$ bits. 
Unfortunately, these tests coincide with 
the existence of the $\mathbb{F}_2$-linear relation \eqref{eqn:relation2} from $21$ to $30$ bits. 
For this, we report the results of the {\it matrix rank} test \cite{Marsaglia1985} and 
the {\it Hamming independence} test \cite{L'Ecuyer:1999:BLC:326147.326156}. 

The {\it matrix rank} test \cite{Marsaglia1985} 
aims to detect linear dependency among blocks of bits. 
The test generates $n$ random $L \times L$ binary matrices, computes 
their ranks, and compares the empirical distribution of these ranks with
their theoretical distribution under $\mathcal{H}_0$ via a chi-square test. 
Each matrix is filled line-by-line by taking 
$\sigma$-bit blocks from $L^2/\sigma$ successive uniform random numbers (assuming that $\sigma$ divides $L$). 
We apply the outputs with $I = \{ 0, 396, 623 \}$ to Test 8 of SmallCrush, 
which is $(n, \tau, \sigma, L)  = (20000, 20, 10, 60)$. 
Table~\ref{table:mt_matrixrank} summarizes the $p$-values. 
As a result, conversion (b) fails such a ``baby" test decisively. 
Note that conversion (b) also fails Tests 57, 59, and 61 of Crush, 
which are the matrix rank tests for the LSBs from $21$ to $30$ bits. 

\begin{table}[t]
\centering
\caption{Summary of $p$-values on the matrix rank test (Test 8 of SmallCrush) for MT19937 with $I = \{ 0, 396, 623 \}$ }
\label{table:mt_matrixrank}
{\footnotesize
\begin{tabular}{l||c|c|c|c|c} \hline
{Conversion} & {1st} & {2nd} & {3rd} & {4th} & {5th} \\ \hline  \hline
{32-bit (a)} & $0.79$ & $0.58$ & $0.79$ & $0.37$ & $0.61$ \\ \hline 
{64-bit (b)} & $< 10^{-300}$ & $< 10^{-300}$ & $< 10^{-300}$ & $< 10^{-300}$ & $< 10^{-300}$ \\ \hline 
\end{tabular}}
\end{table}

Finally, we conduct the {\it Hamming independence} test \cite{L'Ecuyer:1999:BLC:326147.326156}. 
We take $\sigma$ successive bits, say bits $\tau+1$ to $\tau+\sigma$, from each of $2n \lceil L/\sigma \rceil$ successive uniform random numbers, and concatenate these bits to construct $2n$ blocks of $L$ bits.  
Let $X_j$ be the Hamming weight (the number of bits equal to $1$) of the $j$th block 
for $j = 1, \ldots, 2n$. 
Each vector $(X_j, X_{j+1})$ can take $(L + 1) \times (L + 1)$ possible values. 
The test counts the number of occurrences of each possibility among the non-overlapping pairs 
$\{ (X_{2j-1}, X_{2j}), 1 \leq j \leq n \}$. 
These observations are then compared with the expected numbers under 
$\mathcal{H}_0$, via a chi-square test after replacing all classes 
for which the expected number is less than $10$ into a single class. 
For $I = \{ 0, 396, 623\}$, 
we apply (a) and (b) to Test 86 of Crush, 
which is $(n, \tau, \sigma, L) = (10^{8}, 20, 10, 30)$. 
Table~\ref{table:mt_hamming} summarizes the $p$-values. 
Once again, conversion (a) passes but (b) fails the tests with $I = \{ 0, 396, 623\}$. 
Note that each block of $L = 30$ bits is formed 
from the concatenation of $\sigma = 10$ LSBs of 
$\lceil L/\sigma \rceil =3$ successive uniform random numbers 
and hence causes a deviation of the Hamming weight $X_j$ from the binomial distribution $B(L, 1/2)$. 

\begin{table}[t]
\centering
\caption{Summary of $p$-values on the Hamming independence test (Test 86 of Crush) for MT19937 with 
$I = \{ 0, 396, 623 \}$ }
\label{table:mt_hamming}
{\footnotesize
\begin{tabular}{l||c|c|c|c|c} \hline
{Conversion} & {1st} & {2nd} & {3rd} & {4th} & {5th} \\ \hline  \hline
{32-bit (a)} & $0.11$ & $0.22$ & $0.61$ & $0.43$ & $0.85$ \\ \hline 
{64-bit (b)} & $< 10^{-300}$ & $< 10^{-300}$ & $< 10^{-300}$ & $< 10^{-300}$ & $< 10^{-300}$ \\ \hline 
\end{tabular}}
\end{table}

Further, under $I = \{ 0, 396, 623\}$, 
conversion (b) fails the {\it close pairs} test \cite{MR1761436} 
in Test 19 of Crush and Test 22 of BigCrush, which (a) passes. 
Both conversions (a) and (b) fail the {\it linear complexity} test 
in Tests 80 and 81 of BigCrush. 

\begin{remark} \label{remark:decimation}
The 32 MSBs in 
\eqref{eqn:concatenation} are the subsequence obtained
by decimation (i.e., taking a number on every other step) from \eqref{eqn:usual}. 
Hence, conversion (b) with $I = \{ 0, 396, 623 \}$
is essentially the same as conversion (a) with $I = \{ 0, 792, 1246 \}$ in our tests.
However, in general, 
the decimated sequence is not optimized in terms of $k(v)$. 
Thus, when we concatenate the output values, as in \eqref{eqn:concatenation}, 
there is a possibility that the rejection 
occurs in dimensions lower than that of the usual output values, 
that is, half of the dimension in our case. 
This is why we compare 
conversions (a) and (b) with $I = \{ 0, 396, 623 \}$. 
In addition, our result implies that, on 64-bit platforms, 
it is natural to use 64-bit RNGs to avoid such unexpected rejection, 
instead of concatenating 32-bit output values. 

\end{remark}

\begin{remark} \label{remark:res53_tests}
The output from ${\tt genrand\_res53()}$ with $I = \{ 0, 396, 623 \}$ 
in Remark~\ref{remark:mt19937ar} 
fails the birthday spacings tests (Test 12 of Crush and Test 14 of BigCrush), 
overlapping collision test (Test 5 of BigCrush), and linear complexity tests (Tests 80 and 81 of BigCrush), 
but passes the other tests because of the nonexistence of low-weight $\mathbb{F}_2$-linear relations among the LSBs.
\end{remark}



\section{Comparison of MT19937 and other RNGs} \label{sec:comparison}
We now compare MT19937 and other $\mathbb{F}_2$-linear RNGs. 
In the header {\tt <random>} in C++11, 
one can select the 64-bit Mersenne Twister MT19937-64, 
which is an upgraded version of \cite{Nishimura2000} and is available at \url{http://www.math.sci.hiroshima-u.ac.jp/~m-mat/MT/emt64.html}. 
MT19937-64 has $\Delta = 7820$, which is 
much smaller than $\Delta = 13527$ in \eqref{eqn:concatenation}. 
We checked the existence of low-weight $\mathbb{F}_2$-linear relations 
on the $(v \times (k(v)+1))$ bits discussed in Remark~\ref{remark:coding theory}. 
For each $1 \leq v \leq 32$, the minimum weight of 
$\mathbb{F}_2$-linear relations is $\geq 100$, 
so it seems to be safe for the above bits. 
Thus, MT19937-64 is recommended on 64-bit platforms 
when selecting an RNG in the header {\tt <random>}. 

Matsumoto and Saito \cite{SFMT} developed the SIMD-oriented Fast Mersenne Twister SFMT19937.
SFMT19937 has a function to produce 64-bit unsigned integers and 
double-precision floating point numbers in $[0,1)$ by dividing them by $2^{64}$, 
but it has a large $\Delta$ (i.e., $\Delta = 14095$ for 64-bit output) 
and has a similar weakness to MT19937, that is, 
the rejection of the birthday spacings tests with a specific lag set $I$. 
(See the online appendix of \cite{Harase:2018:IME:3175005.3159444} for details.) 
Matsumoto and Saito \cite{MR2743921} 
developed the dSFMT generator, which is a 52-bit generator specialized for generating
double-precision floating-point numbers based on the IEEE 754 format 
using a union (i.e., the 12 bits for the sign and exponent are fixed, and the 52 bits of the significand are taken from the generator).  
This generator has a small $\Delta$ (i.e., $\Delta =2616$ for 52-bit output).

Harase and Kimoto \cite{Harase:2018:IME:3175005.3159444} developed 
64-bit maximally equidistributed 
$\mathbb{F}_2$-linear generators with Mersenne prime periods from $2^{521}-1$ to $2^{44497}-1$, 
named 64-bit MELGs. 
For the 64-bit MELGs, 
the minimum weight of the $\mathbb{F}_2$-linear relations among the $(v \times (k(v)+1))$ bits 
is approximately (or slightly smaller than) $p/2$ for each $1 \leq v \leq32$, e.g., 
MELG19937-64 (with period $2^{19937}-1$)
has a minimum weight greater than $9500$. 
For the 32 LSBs,
the existence of low-weight $\mathbb{F}_2$-linear relations (e.g., number of terms $< 20$) is also avoided. 
Harase and Kimoto implemented a function to generate double-precision floating-point numbers in $[0,1)$ using a similar union to that in the dSFMT generator. 

L'Ecuyer \cite{MR1620231} developed the $64$-bit combined linear feedback shift register generator LFSR258 with a small state space. In the above setting, this generator also avoids the low-weight $\mathbb{F}_2$-linear relations for the 32 MSBs and 32 LSBs, respectively.

\section{Concluding remarks} \label{sec:conclusion}
For 32-bit RNGs, the concatenation of outputs naturally appears in the conversion to
53-bit double-precision floating-point numbers. 
However, the RNGs are not usually optimized for this purpose, 
so the concatenation might degrade 
the dimension of equidistribution $k(v)$ compared with the usual 32-bit use. 
In this paper, we have analyzed the concatenation of outputs 
generated by the 32-bit MT19937. 
As pointed out in \cite{MR3191969}, MT19937 has low-weight 
$\mathbb{F}_2$-linear relations in just one dimension 
beyond $k(v)$. 
When concatenating the outputs, 
the above $\mathbb{F}_2$-linear relations are folded, 
which destroys the uniformity, causing the $k(v)$ values to decrease. 
In addition, we showed the existence of a low-weight $\mathbb{F}_2$-linear relation 
among the LSBs of MT19937. 
We observed dependency among bits in several statistical tests. 

For most real-life simulations, these properties are not likely to affect the results, 
because the simulation models are not synchronized with the existence of 
low-weight $\mathbb{F}_2$-linear relations (unless we are very unlucky). 
However, such a structure should be avoided whenever possible when designing a new RNG \cite{MR1478040,MR3246605,L'Ecuyer:2004:DRN:961292.961302}.

In summary, we can say that 
RNGs with recurrence based on 64-bit words are preferable on 64-bit platforms. 
In the header {\tt <random>} of the C++11 STL, 
MT19937-64, a 64-bit version of MT, has been implemented, 
and this 64-bit RNG is recommended 
when converting into double-precision floating-point numbers in $[0,1)$. 
As alternatives, we note that 64-bit maximally equidistributed 
$\mathbb{F}_2$-linear RNGs are now available \cite{Harase:2018:IME:3175005.3159444,MR1620231}.

\subsection*{Acknowledgments}
The author would like to thank the reviewers for their valuable comments and suggestions. 
This work was partially supported by JSPS Grant-in-Aid for Early-Career Scientists $\#$18K18016, 
for Young Scientists (B) $\#$26730015, for Scientific Research (B) $\#$26310211, 
for Challenging Exploratory Research $\#$15K13460, and by Research Promotion Program 
for Acquiring Grants-in-Aid for Scientific Research at Ritsumeikan University.
This work was also supported by JST CREST. 




\bibliographystyle{model1b-num-names.bst}
\bibliography{harase-mt}   
\end{document}